\documentstyle[12pt]{amsart}

\def\opn#1#2{\def#1{\operatorname{#2}}}
\opn\reg{reg}
\newtheorem{Theorem}{Theorem}[section]
\newtheorem{Lemma}[Theorem]{Lemma}

\newtheorem{Proposition}[Theorem]{Proposition}

\begin{document}

\title{on the Betti numbers of shifted complexes of stable simplicial complexes}
\author{Zhongming Tang and Guifen Zhuang}

\begin{abstract}
Let $\Delta$ be a stable simplicial complex on $n$ vertexes. Over
an arbitrary base field $K$, the symmetric algebraic shifted
complex $\Delta^s$ of $\Delta$ is defined. It is proved that the
Betti numbers of the Stanley-Reisner ideals in the polynomial ring
$K[x_1,x_2,\ldots,x_n]$ of the symmetric algebraic shifted,
exterior algebraic shifted and combinatorial shifted complexes of
$\Delta$ are equal.
\end{abstract}
\thanks{
\indent Supported by the National Natural Science Foundation of
China}
\date{}
\maketitle

\section*{Introduction}
Throughout the paper, let $\Delta$ be a simplicial complex on the
vertex set $[n]=\{1,2,\ldots,n\}$. Then, by definition, $\Delta$
is a set of some subsets of $[n]$ such that $\{i\}\in\Delta$,
$i=1,2,\ldots,n$, and, for any $\sigma\in\Delta$, if
$\tau\subseteq\sigma$ then $\tau\in\Delta$. Set
$\nabla=\{\sigma\subseteq[n] : \sigma\not\in\Delta\}$. Let
$S=K[x_1,x_2,\ldots,x_n]$ be the polynomial ring over an infinite
field $K$. For any $\sigma=\{i_1,i_2,\ldots,i_r\}\subseteq[n]$
with $i_1<i_2<\ldots<i_r$, we set $x_{\sigma}=x_{i_1}x_{i_2}\cdots
x_{i_r}$. Then $\Delta$ determines an ideal $I_{\Delta}$, the
Stanley-Reisner ideal of $\Delta$, of $S$, which is the monomial
ideal generated by all monomials $x_{\sigma}$ such that
$\sigma\in\nabla$. Conversely, for any monomial ideal $L$ of $S$
which is generated by some squarefree monomials, there is a unique
simplicial complex $\Delta_1$ on the vertex set $[n]$ such that
$L=I_{\Delta_1}$. We refer the reader to \cite{BH} and \cite{S}
for the detailed information about these notions.

For any graded ideal $I$ of $S$, the graded Betti number
$\beta_{i,j}(I)$ is defined by the minimal graded free
$S$-resolution of $I$:
$$
\cdots\rightarrow\oplus_jS(-j)^{\beta_{i,j}(I)}\rightarrow\cdots
\rightarrow\oplus_jS(-j)^{\beta_{1,j}(I)}\rightarrow\oplus_jS(-j)^{\beta_{0,j}(I)}
\rightarrow I\rightarrow 0.
$$
Set $\beta_i(I)=\sum_j\beta_{i,j}(I)$. The regularity, $\reg(I)$,
of $I$ is defined as $\reg(I)=\max\{j-i : \beta_{i,j}(I)\not=0\}$.

For the simplicial complex $\Delta$, there are three shifted
complexes $\Delta^s$, $\Delta^e$ and $\Delta^c$, which are called
the symmetric algebraic shifted complex, exterior algebraic
shifted complex and combinatorial shifted complex of $\Delta$
respectively, cf. \cite{AHH1}. Shifted complexes are determined by
some monomial ideals of a polynomial ring or an exterior algebra
over a field $K$. The combinatorial shifting and the exterior
algebraic shifting are defined over any base field $K$, while the
symmetric algebraic shifting is only defined when
$\mbox{char}(K)=0$, cf. \cite{AHH1}. In the present paper, we
firstly show that if $I$  is a (squarefree) stable ideal of $S$,
then $\mbox{Gin}(I)$ is also stable. Then it makes sense to define
the symmetric algebraic shifting of stable complexes over an
arbitrary base field $K$. One of the important problems in the
study of simplicial complexes is the behaviour of the graded Betti
numbers under shifting. In \cite{AHH1}, Aramova, Herzog and Hibi
conjectured that
$\beta_{i,j}(I_{\Delta^s})\leq\beta_{i,j}(I_{\Delta^e})\leq\beta_{i,j}(I_{\Delta^c})$,
although their symmetric algebraic shifting is only defined in
characteristic $0$. In this paper, we will show that, over any
base field, when $\Delta$ is stable,
$\beta_{i,j}(I_{\Delta^s})=\beta_{i,j}(I_{\Delta^e})=\beta_{i,j}(I_{\Delta^c})=\beta_{i,j}(I_{\Delta})$,
for all $i$ and $j$.

\medskip
\section{Betti numbers under symmetric algebraic shifting}
A monomial $x_1^{a_1}x_2^{a_2}\cdots x_n^{a_n}\in S$ is said to be
squarefree if $a_i\leq 1$, $i=1,2,\ldots,n$. For a monomial $u\in
S$, we denote $\max\{i : x_i\mid u\}$ by $m(u)$, and, for any
$\sigma=\{i_1,i_2,\ldots,i_r\}\subseteq[n]$, we also denote
$\max\{i_s : s=1,2,\ldots,r\}$ by $m(\sigma)$. Let $I$ be a
monomial ideal of $S$, the minimal system of monomial generators
of $I$ will be denoted by $G(I)$. For any $j$, we write $I_j$ (
$G(I)_j$ ) for the set of monomials of degree $j$ belonging to $I$
( $G(I)$ )and $I_{\langle j\rangle}$ for the ideal generated by
$I_j$. A monomial ideal $I$ of $S$ is said to be squarefree if it
is generated by squarefree monomials. Then, it is clear that
$I_{\Delta}$ is a squarefree monomial ideal.

A (squarefree) monomial ideal is said to be (squarefree) stable if
its set $B$ of (squarefree) monomials is (squarefree) stable,
i.e., $x_i(u/x_{m(u)})\in B$ for all $u\in B$ and all $i<m(u)$ (
such that $x_i \not| u$ ), further, if $x_i(u/x_j)\in B$ for all
$u\in B$ and all $i<j$ such that $x_j\mid u$ ( and $x_i \not| u$ )
then we say that the ideal is strongly (squarefree) stable. Let
$I$ be a monomial ideal of $S$. If $I$ is stable, it is shown in
\cite{EK} that
$$
\beta_{i,i+j}(I)=\sum_{u\in G(I)_j}\left(
\begin{array}{c}
m(u)-1\\
i
\end{array}
\right)
$$
and
$$
\reg(I)=\max\{\deg(u):u\in G(I)\}.
$$
If $I$ is squarefree stable, then it is given in \cite{AHH3} that
$$
\beta_{i,i+j}(I)=\sum_{u\in G(I)_j}\left(
\begin{array}{c}
m(u)-j\\
i
\end{array}
\right)
$$
and
$$
\reg(I)=\max\{\deg(u):u\in G(I)\}.
$$

We say that the simplicial complex $\Delta$ is stable if
$I_{\Delta}$ is squarefree stable. Then, it is easy to show the
following

\begin{Lemma} The following are equivalent
\begin{itemize}
\item[(1)] $\Delta$ is stable.
\item[(2)] For any $\sigma\in\nabla$, if $1\leq i<m(\sigma)$ and
$i\not\in\sigma$, then
$(\sigma\setminus\{m(\sigma)\})\cup\{i\}\in\nabla$.
\item[(3)] For any $\sigma\in\Delta$, if $i\in\sigma$ and
$m(\sigma)<j\leq n$, then
$(\sigma\setminus\{i\})\cup\{j\}\in\Delta$.
\end{itemize}
\end{Lemma}

For a graded ideal $I$ of $S$, let $\mbox{Gin}(I)$ be the generic
initial ideal of $I$ with respect to the reverse lexicographic
term order on $S$ induced by $x_1>x_2>\ldots>x_n$. Then
$\mbox{Gin}(I)$ is Borel-fixed and, when $\mbox{char}(K)=0$, it is
strongly stable, cf. \cite[Chapter 15]{E}. For a general field
$K$, we have the following

\begin{Proposition}
\label{ginstable} If $I$  is a (squarefree) stable ideal of $S$,
then ${\em \mbox{Gin}}(I)$ is also stable.
\end{Proposition}

\begin{pf}
Let $j\geq 1$ be such that $I_j\not=0$. As $I$ is stable or
squarefree stable, it follows that $I$ is componentwise linear,
hence $I_{\langle j\rangle}$ is $j$-linear. Note that
$\mbox{reg}(I_{\langle j\rangle})=j$. Then
$\mbox{reg}(\mbox{Gin}(I_{\langle
j\rangle}))=\mbox{reg}(I_{\langle j\rangle})=j$, thus
$\mbox{Gin}(I_{\langle j\rangle})$ is generated in degree $j$. By
the construction of generic initial ideals, we see that
$\mbox{Gin}(I)_j=\mbox{Gin}(I_{\langle j\rangle})_j$, then it
turns out that $\mbox{Gin}(I)_{\langle
j\rangle}=\mbox{Gin}(I_{\langle j\rangle})$. Since
$\mbox{Gin}(I_{\langle j\rangle})$ is Borel-fixed and generated in
degree $j$ and $\mbox{reg}(\mbox{Gin}(I_{\langle j\rangle}))=j$,
it follows from \cite[Proposition 10]{ERT} that
$\mbox{Gin}(I_{\langle j\rangle})$ is stable, i.e.,
$\mbox{Gin}(I)_{\langle j\rangle}$ is stable. Hence
$\mbox{Gin}(I)$ is stable.
\end{pf}

There is an operator $\sigma$ which transfers any monomial ideal
into a squarefree monomial ideal. For any $u=x_{i_1}x_{i_2}\cdots
x_{i_j}\cdots x_{i_d}\in S$ where $i_1\leq i_2\leq\cdots\leq
i_j\leq\cdots\leq i_d$, we set $u^{\sigma}=x_{i_1}x_{i_2+1}\cdots
x_{i_j+(j-1)}\cdots x_{i_d+(d-1)}$. Let $I$ be a monomial ideal of
$S$ with $G(I)=\{u_1,u_2,\ldots,u_r\}$. Then we set $I^{\sigma}$
for the squarefree monomial ideal generated by
$u^{\sigma}_1,u_2^{\sigma},\ldots,u_r^{\sigma}$. Thus $I^{\sigma}$
is a squarefree ideal of $K[x_1,x_2,\ldots,x_m]$, where
$m=\max\{m(u)+\deg(u)-1 : u\in G(I)\}$.

Suppose that $\Delta$ is  stable. Then, by \ref{ginstable},
$\mbox{Gin}(I_{\Delta})$ is also stable. Using the same arguments
as in the proof of \cite[Lemma 1.1]{AHH1}, we see that
$m(u)+\deg(u)-1\leq n$ for all $u\in G(\mbox{Gin}(I_{\Delta}))$,
hence, $(\mbox{Gin}(I_{\Delta}))^{\sigma}$ is also an ideal of
$S$. Then we define the symmetric algebraic shifted complex
$\Delta^s$ of $\Delta$ by
$$
I_{\Delta^s}=(\mbox{Gin}(I_{\Delta}))^{\sigma}.
$$

We are going to show that
$\beta_{i,j}(I_{\Delta})=\beta_{i,j}(I_{\Delta^s})$ if $\Delta$ is
stable. When $\mbox{char}(K)=0$, this result can be easily proved.
But, for a base field $K$ with arbitrary characteristic, we need
some new results as we do not know whether
$\mbox{Gin}(I_{\Delta})$ is strongly stable.

\begin{Proposition}
\label{ginbetti} If $I$ is a (squarefree) stable ideal of $S$,
then
$$
\beta_{i,j}(I)=\beta_{i,j}({\em \mbox{Gin}}(I)) \mbox{  for all
$i$ and $j$ },
$$
\end{Proposition}

\begin{pf}
By \ref{ginstable},  $\mbox{Gin}(I)$ is also stable, thus $I$ and
$\mbox{Gin}(I)$ are componentwise linear. Then, by
\cite[Proposition 1.3]{HH}, for any $i$ and $j$,
\begin{eqnarray*}
\beta_{i,i+j}(I)&=& \beta_i(I_{\langle j\rangle})-\beta_i({\frak m}I_{\langle j-1\rangle}),\\
\beta_{i,i+j}( \mbox{Gin}(I))&=& \beta_i( \mbox{Gin}(I)_{\langle
j\rangle})-
\beta_i({\frak m} \mbox{Gin}(I)_{\langle j-1\rangle})\\
&=&\beta_i( \mbox{Gin}(I_{\langle j\rangle}))-\beta_i( {\frak
m}\mbox{Gin}(I_{\langle j-1\rangle})),
\end{eqnarray*}
where ${\frak m}$ is the ideal $(x_1,x_2,\ldots,x_n)$ of $S$.
Since ${\frak m}I_{\langle j-1\rangle}$ is $j$-linear, it follows
that $\mbox{reg}({\frak m}I_{\langle j-1\rangle})=j$, hence
$\mbox{reg}(\mbox{Gin}({\frak m}I_{\langle j-1\rangle}))=j$. Then
$\mbox{Gin}({\frak m}I_{\langle j-1\rangle})$ is also generated in
degree $j$. On the other hand, as $\mbox{Gin}(I_{\langle
j-1\rangle})$ is generated in degree $j-1$ (refer to the proof of
\ref{ginstable}), we see that ${\frak m}\mbox{Gin}(I_{\langle
j-1\rangle})$ is also generated in degree $j$. But
$$
\mbox{Gin}({\frak m}I_{\langle
j-1\rangle})_j=\mbox{Gin}(I_{\langle j-1\rangle})_j=({\frak
m}\mbox{Gin}(I_{\langle j-1\rangle}))_j,
$$
hence ${\frak m}\mbox{Gin}(I_{\langle
j-1\rangle})=\mbox{Gin}({\frak m}I_{\langle j-1\rangle})$. To end
the proof, it is enough to show that, if $I$ is (squarefree)
stable and generated in degree $d$, then
$\beta_i(I)=\beta_i(\mbox{Gin}(I))$. By the assumptions, $I$ and
$\mbox{Gin}(I)$ are $d$-linear and their Hilbert functions have
the following forms:
\begin{eqnarray*}
H_I(t)&=&(\sum_i(-1)^i\beta_{i,i+d}(I)t^i)\frac{t^d}{(1-t)^n},\\
H_{{\rm
Gin}(I)}(t)&=&(\sum_i(-1)^i\beta_{i,i+d}({\mbox{Gin}(I)})t^i)\frac{t^d}{(1-t)^n}.
\end{eqnarray*}
But $I$ and $\mbox{Gin}(I)$ have the same Hilbert function, hence
$\beta_i(I)=\beta_i(\mbox{Gin}(I))$, as required.
\end{pf}

The following result is proved in \cite{AHH1} when $I$ is strongly
stable.

\begin{Lemma}
\label{sigmabetti} If $I$ is stable and Borel-fixed, then
$\beta_{i,j}(I)=\beta_{i,j}(I^{\sigma})$, for all $i$ and $j$.
\end{Lemma}

\begin{pf}
Let $G(I)=\{u_1,u_2,\ldots,u_r\}$. We will show that $I^{\sigma}$
is squarefree stable with
$G(I^{\sigma})=\{u_1^{\sigma},u_2^{\sigma},\ldots,u_r^{\sigma}\}$.
Then the result follows from the formulas of Betti numbers of
stable and squarefree stable ideals.

Let $u=x_{i_1}^{r_1}x_{i_2}^{r_2}\cdots x_{i_s}^{r_s}\in I$ with
$i_1<i_2<\ldots<i_s$ and $r_j\geq 1$, $j=1,2,\ldots,s$. We firstly
note that, since $I$ is $p$-Borel with $p=\mbox{char}(K)$ and
$k\leq_pk$ for any positive integer $k$, cf. \cite[Chapter 15]{E},
it follows that $(x_k/x_{i_j})^{r_j}u\in I$ for all $1\leq k<i_j$.
Secondly, let $v=x_{j_1}^{t_1}x_{j_2}^{t_2}\cdots x_{j_l}^{t_l}\in
I$ with $j_1<j_2<\ldots<j_l$ and $t_k\geq 1$, $k=1,2,\ldots,l$.
Suppose that $v^{\sigma}\mid u^{\sigma}$ and to see the
relationship between $u$ and $v$. For any $r\in[n]$ and
$J\subseteq[n]\cup\{0\}$, set $x_{r+J}=\prod_{k\in J}x_{r+k}$ and
$[0,r-1]=\{0,1,\ldots,r-1\}$. Then
$$
u^{\sigma}=x_{i_1+[0,r_1-1]}\cdot x_{i_2+r_1+[0,r_2-1]}\cdots
x_{i_s+(r_1+\ldots+r_{s-1})+[0,r_s-1]}.
$$
We notice that the difference between indices of any two $x_w$ in
different groups\\
 $x_{i_k+(r_1+\ldots+r_{k-1})+[0,r_k-1]}$ is at
least 2. Thus, from
$$
v^{\sigma}=x_{j_1+[0,t_1-1]}\cdot x_{j_2+t_1+[0,t_2-1]}\cdots
x_{j_l+(t_1+\ldots+t_{l-1})+[0,t_l-1]}
$$
and $v^{\sigma}\mid u^{\sigma}$, we see that $l\leq s$ and there
exist $w_k\in[s]$, $k=1,2,\ldots,l$, such that
$$
x_{j_k+(t_1+\ldots+t_{k-1})+[0,t_k-1]}\mid
x_{i_{w_k}+(r_1+\ldots+r_{w_k-1})+[0,r_{w_k}-1]}, k=1,2,\ldots,l.
$$
Hence $j_k=i_{w_k}+m_k$ with $m_k\geq 0$ and $t_k\leq r_{w_k}$,
$k=1,2,\ldots,l$. It turns out that
$x_{j_1-m_1}^{t_1}x_{j_2-m_2}^{t_2}\cdots x_{j_l-m_l}^{t_l}\mid u$
for some $m_k\geq 0$, $k=1,2,\ldots,l$.

Now we show that
$G(I^{\sigma})=\{u_1^{\sigma},u_2^{\sigma},\ldots,u_r^{\sigma}\}$.
Suppose that $u_k^{\sigma}\not\in G(I^{\sigma})$. Then there
exists some $l\not=k$ such that $u_l^{\sigma}\mid u_k^{\sigma}$.
Let $u_l=x_{i_1}^{r_1}x_{i_2}^{r_2}\cdots x_{i_s}^{r_s}$ with
$i_1<i_2<\ldots<i_s$ and $r_j\geq 1$, $j=1,2,\ldots,s$. Then
$x_{i_1-m_1}^{r_1}x_{i_2-m_2}^{r_2}\cdots x_{i_s-m_s}^{r_s}\mid
u_k$ for some $m_t\geq 0$, $t=1,2,\ldots,s$. Note that
$x_{i_1-m_1}^{r_1}x_{i_2-m_2}^{r_2}\cdots
x_{i_s-m_s}^{r_s}\not=u_k$ by $u_l^{\sigma}\mid u_k^{\sigma}$.
Since $x_{i_1}^{r_1}x_{i_2}^{r_2}\cdots x_{i_s}^{r_s}\in I$ and
$I$ is Borel-fixed, as noted above, we see that
$x_{i_1-m_1}^{r_1}x_{i_2-m_2}^{r_2}\cdots x_{i_s-m_s}^{r_s}\in I$.
This contradicts that $u_k\in G(I)$. Thus
$G(I^{\sigma})=\{u_1^{\sigma},u_2^{\sigma},\ldots,u_r^{\sigma}\}$.

Finally, we show that $I^{\sigma}$ is squarefree stable. Let
$u=x_{i_1}^{r_1}x_{i_2}^{r_2}\cdots x_{i_s}^{r_s}\in G(I)$ with
$i_1<i_2<\ldots<i_s$ and $r_j\geq 1$, $j=1,2,\ldots,s$. Suppose
that $l<m(u^{\sigma})=i_s+(r_1+\ldots+r_s)-1$ and $x_l \not|
u^{\sigma}$. Let
$i_{k-1}+(r_1+\ldots+r_{k-1})-1<l<i_k+(r_1+\ldots+r_{k-1})$ and
set $v=x_{i_1}^{r_1}\cdots
x_{i_{k-1}}^{r_{k-1}}x_{l-(r_1+\ldots+r_{k-1})}x_{i_k-1}^{r_k}\cdots
x_{i_{s-1}-1}^{r_{s-1}} x_{i_s-1}^{r_s-1}$. Note that $v\in I$ as
above and $v^{\sigma}=(x_l/x_{m(u^{\sigma})})u^{\sigma}$. Then
$v=u_wv^*$ for some $u_w\in G(I)$ and a monomial $v^*\in S$, as
$I$ is stable, we may assume that $m(u_w)<j$ for all $j$ such that
$x_j\mid v^*$. Thus we have that $u_w^{\sigma}\mid v^{\sigma}$,
hence $v^{\sigma}\in I^{\sigma}$, as required.
\end{pf}

Now we prove the following

\begin{Theorem}
\label{sbetti} If $\Delta$ is stable, then, for all $i$ and $j$,
$$
\beta_{i,j}(I_{\Delta})=\beta_{i,j}(I_{\Delta^s}).
$$
\end{Theorem}

\begin{pf}
Since $I_{\Delta}$ is squarefree stable, it follows from
\ref{ginbetti} that
$\beta_{i,j}(I_{\Delta})=\beta_{i,j}(\mbox{Gin}(I_{\Delta}))$. But
$\mbox{Gin}(I_{\Delta})$ is also stable by \ref{ginstable}, and
Borel-fixed. Then the theorem follows from \ref{sigmabetti} .
\end{pf}

For exterior algebraic shifting, a similar result can been easily
obtained. Let $V$ be an $n$-dimensional $K$-vector space with
basis $e_1,e_2,\ldots,e_n$. Then the exterior algebra
$E=\bigwedge^{\cdot}(V)$ is a finite dimensional graded
$K$-algebra and its $r$-th graded component $\bigwedge^r(V)$ has
the $K$-basis $e_{i_1}\wedge e_{i_2}\wedge\cdots\wedge e_{i_r}$
with $i_1<i_2<\ldots<i_r$. For any
$\sigma=\{i_1,i_2,\ldots,i_r\}\subseteq[n]$ with
$i_1<i_2<\ldots<i_r$, we set $e_{\sigma}=e_{i_1}\wedge
e_{i_2}\wedge\cdots\wedge e_{i_r}$. The stable and strongly stable
ideals of $E$ ( which are automatically squarefree ) can be
defined as in $S$. Let $I$ be a graded ideal of $E$. The generic
initial ideal, $\mbox{Gin}(I)$, of $I$ in $E$ is also defined
similarly as in $S$. Then  $\mbox{Gin}(I)$ is strongly stable for
any infinite field $K$, cf. \cite{AHH4}.

For the simplicial complex $\Delta$, we set $J_{\Delta}$, the
Stanley-Reisner ideal of $\Delta$ in $E$, as the monomial ideal of
$E$ generated by all monomials $e_{\sigma}$ such that
$\sigma\in\nabla$. The exterior algebraic shifted complex
$\Delta^e$ of $\Delta$ is defined by
$$
J_{\Delta^e}=\mbox{Gin}(J_{\Delta}).
$$

Suppose that $\Delta$ is stable. Then $J_{\Delta}$ is also
componentwise linear. Thus, by \cite[Corollary 2.2]{AHH2}, we have
the following

\begin{Proposition}
\label{ebetti} If $\Delta$ is stable, then, for all $i$ and $j$,
$$
\beta_{i,j}(I_{\Delta})=\beta_{i,j}(I_{\Delta^e}).
$$
\end{Proposition}

\medskip
\section{Betti numbers under combinatorial shifting}
Let $1\leq k<l\leq n$, for any $\sigma\in\nabla$, we define
$$
S_{kl}(\sigma)=\left\{
\begin{array}{ll}
(\sigma\setminus\{l\})\cup\{k\}&\mbox{if $l\in\sigma$,
$k\not\in\sigma$ and $(\sigma\setminus\{l\})\cup\{k\}\not\in\nabla$},\\
\sigma&\mbox{otherwise}.
\end{array}
\right.
$$
We summarize  some properties of $S_{kl}$ as follows:
\begin{itemize}
\item[(1)] For any $\sigma\in\nabla$,
$S_{kl}(\sigma)\not\in\nabla$ if $S_{kl}(\sigma)\not=\sigma$.
\item[(2)] Let $\sigma,\tau\in\nabla$, if $\sigma\not=\tau$ then $S_{kl}(\sigma)\not= S_{kl}(\tau)$.
\item[(3)] If $\Delta$ is stable, then $m(S_{kl}(\sigma))=m(\sigma)$.
\item[(4)] Let
$G(I_{\Delta})=\{x_{\sigma_1},x_{\sigma_2},\ldots,x_{\sigma_r}\}$.
Then $\sigma_i\not=S_{kl}(\sigma_j)$ for any $i\not=j$, and, if
$S_{kl}(\sigma_i)\subseteq S_{kl}(\sigma_j)$ for $i\not=j$, then
$S_{kl}(\sigma_i)\not=\sigma_i$ and $S_{kl}(\sigma_j)=\sigma_j$.
\end{itemize}

$\mbox{Shift}_{kl}(\Delta)$ is the simplicial complex defined by
$I_{{\rm Shift}_{kl}(\Delta)}$ which is the squarefree monomial
ideal of $S$ generated by all monomials $x_{S_{kl}(\sigma)}$ with
$\sigma\in\nabla$. A combinatorial shifted complex $\Delta^c$ of
$\Delta$ is a shifted complex and
$$
\Delta^c=\mbox{Shift}_{k_rl_r}(\mbox{Shift}_{k_{r-1}l_{r-1}}(\cdots(\mbox{Shift}_{k_1l_1}(\Delta))\cdots)),
$$
for some $r\geq 1$ and $k_i<l_i$, $i=1,2,\ldots,r$.

We firstly discuss the properties of $\mbox{Shift}_{kl}(\Delta)$
and, then, get the properties of $\Delta^c$ by induction.

\begin{Lemma}
\label{cstable} If $\Delta$ is stable, then ${\em
\mbox{Shift}}_{kl}(\Delta)$ is also stable.
\end{Lemma}

\begin{pf}
It is enough to show that, for any $\sigma\in\nabla$, if
$i<m(S_{kl}(\sigma))$ and $i\not\in S_{kl}(\sigma)$, then
$(S_{kl}(\sigma)\setminus\{m(S_{kl}(\sigma))\})\cup\{i\}=S_{kl}(\sigma')$
for some $\sigma'\in\nabla$. Set
$\sigma_i=(S_{kl}(\sigma)\setminus\{m(S_{kl}(\sigma))\})\cup\{i\}$.
Consider two cases.

Case I)  $S_{kl}(\sigma)=\sigma$. Then
$\sigma_i=(\sigma\setminus\{m(\sigma)\})\cup\{i\}\in\nabla$ by the
stability of $\Delta$. We will show that
$\sigma_i=S_{kl}(\sigma_i)$. On the contrary, suppose that
$S_{kl}(\sigma_i)\not=\sigma_i$. Then
$S_{kl}(\sigma_i)=(\sigma_i\setminus\{l\})\cup\{k\}$,
$l\in\sigma_i$, $k\not\in\sigma_i$ and
$(\sigma_i\setminus\{l\})\cup\{k\}\not\in\nabla$. Note that
$i\not=k$ and
$(\sigma_i\setminus\{l\})\cup\{k\}=(((\sigma\setminus\{m(\sigma)\})\cup\{i\})\setminus\{l\})\cup\{k\}$.
Consider two possibilities: $i=l$ and $i\not=l$. When $i=l$, as
$\Delta$ is stable, $k<m(\sigma)$ and $k\not\in\sigma$, we have
that
$$
(\sigma_i\setminus\{l\})\cup\{k\}=(\sigma\setminus\{m(\sigma)\})\cup\{k\}\in
\nabla,
$$
a contradiction. When $i\not=l$, as $l\in\sigma_i$,
$k\not\in\sigma_i$ and $k<l$, we see that $l\in \sigma$ and
$k\not\in\sigma$. But $S_{kl}(\sigma)=\sigma$, this implies that
$(\sigma\setminus\{l\})\cup\{k\}\in\nabla$. By $l\in\sigma_i$
again, we have that $l\not=m(\sigma)$, hence
$m((\sigma\setminus\{l\})\cup\{k\})=m(\sigma)$ and
$$
(\sigma_i\setminus\{l\})\cup\{k\}=(((\sigma\setminus\{l\})\cup\{k\})
\setminus\{m((\sigma\setminus\{l\})\cup\{k\})\})\cup\{i\}.
$$
Since $i<m(\sigma)$ and $i\not\in(\sigma\setminus\{l\})\cup\{k\}$,
it follows from the stability of $\Delta$ that the right-hand side
of the above equality is in $\nabla$, i.e.,
$(\sigma_i\setminus\{l\})\cup\{k\}\in\nabla$, a contradiction
again.

Case II) $S_{kl}(\sigma)\not=\sigma$. Then
$S_{kl}(\sigma)=(\sigma\setminus\{l\})\cup\{k\}$, $l\in\sigma$,
$k\not\in\sigma$ and
$(\sigma\setminus\{l\})\cup\{k\}\not\in\nabla$. Thus, by the
stability of $\Delta$, we see immediately that $l\not=m(\sigma)$,
hence $m(S_{kl}(\sigma))=m(\sigma)$ and
$\sigma_i=(((\sigma\setminus\{l\})\cup\{k\})\setminus\{m(\sigma)\})\cup\{i\}$.
Note that $i\not=k$ by $i\not\in S_{kl}(\sigma)$. Consider two
possibilities: $i=l$ and $i\not=l$ again. When $i=l$, as
$k<m(\sigma)$ and $k\not\in\sigma$, we have that
$\sigma_i=(\sigma\setminus\{m(\sigma)\})\cup\{k\}\in\nabla$ by the
stability of $\Delta$, and $S_{kl}(\sigma_i)=\sigma_i$, as
required. When $i\not=l$, we have that
$\sigma_i=(((\sigma\setminus\{m(\sigma)\})\cup\{i\})\setminus\{l\})\cup\{k\}$.
As $i<m(\sigma)$ and $i\not\in\sigma$, we see that
$(\sigma\setminus\{m(\sigma)\})\cup\{i\}\in\nabla$ by the
stability of $\Delta$. Since
$l\in(\sigma\setminus\{m(\sigma)\})\cup\{i\}$ and
$k\not\in(\sigma\setminus\{m(\sigma)\})\cup\{i\}$ by
$k\not\in\sigma$ and $i\not=k$, it follows that
$\sigma_i=S_{kl}(\sigma_i)$ if $\sigma_i\in\nabla$ and
$\sigma_i=S_{kl}((\sigma\setminus\{m(\sigma)\})\cup\{i\})$ if
$\sigma_i\not\in\nabla$, as required.
\end{pf}

Suppose that $\Delta$ is stable and $G(I_{\Delta})=\{x_{\sigma_1},
x_{\sigma_2},\ldots, x_{\sigma_r}\}$. For any $\sigma\in\nabla$,
there is an $i\in[r]$ such that $\sigma=\sigma_i\cup\tau$ with
$\tau\subseteq[n]$, where we understand always that
$\sigma\cap\tau=\emptyset$. Let $\tau_1,\tau_2\subseteq[n]$, we
write $\tau_1<\tau_2$ if $s<t$ for all $s\in\tau_1$ and
$t\in\tau_2$. Note that, as $\Delta$ is stable, we can find some
$j\in[r]$ such that $\sigma=\sigma_j\cup\tau'$ where
$\tau'\subseteq[n]$ and $\sigma_j<\tau'$, and such a
representation is unique.

Let $i\in[n]$. If $S_{kl}(\sigma_i)\supseteq S_{kl}(\sigma_j)$ for
some $j\not=i$, then $S_{kl}(\sigma_i)=\sigma_i$ and
$S_{kl}(\sigma_j)\not=\sigma_j$, hence
$\sigma_i=S_{kl}(\sigma_j)\cup\tau$ with $\tau\subseteq[n]$ and
$\tau\not=\emptyset$, and conversely, this equality implies that
$S_{kl}(\sigma_j)\not=\sigma_j$, $k\in\sigma_i$, hence
$S_{kl}(\sigma_i)=\sigma_i$ and $S_{kl}(\sigma_i)\supseteq
S_{kl}(\sigma_j)$. Furthermore, $\sigma_j\cup\tau$ is uniquely
defined by $\sigma_i$ since, if
$\sigma_i=S_{kl}(\sigma_{j'})\cup\tau'$ with $j'\in[r]$ and
$\tau'\subseteq[n]$ then $l\in\sigma_j,\sigma_{j'}$,
$k\not\in\sigma_j,\sigma_{j'}$ and
$((\sigma_j\setminus\{l\})\cup\{k\})\cup\tau=((\sigma_{j'}\setminus\{l\})\cup\{k\})\cup\tau'$,
hence $\sigma_j\cup\tau=\sigma_{j'}\cup\tau'$ recalling that, by
our convention,
$S_{kl}(\sigma_j)\cap\tau=S_{kl}(\sigma_{j'})\cap\tau'=\emptyset$.
As noted above, we may assume that $\sigma_j<\tau$. Then it makes
sense to give the following definition.

Let $i\in[r]$, we define $\sigma_i^*$ as follows
$$
\sigma_i^*=\left\{
\begin{array}{ll}
\sigma_s\cup\tau&\mbox{ if $\sigma_i=S_{kl}(\sigma_s)\cup\tau$ for
some $s\not=i$}\\
S_{kl}(\sigma_i)&\mbox{ otherwise}
\end{array}
\right.
$$
Thus there are three possible values for $\sigma_i^*$. If
$S_{kl}(\sigma_i)\not=\sigma_i$ then
$\sigma_i^*=S_{kl}(\sigma_i)\not\in\nabla$; if
$S_{kl}(\sigma_i)=\sigma_i$ then $\sigma_i^*=\sigma_s\cup\tau$ or
$\sigma_i$ depending on whether
$\sigma_i=S_{kl}(\sigma_s)\cup\tau$ for some $s\not=i$ or not.
Note that $m(\sigma_i^*)=m(\sigma_i)$, $i=1,2,\ldots,r$.

\begin{Lemma}
\label{genarator} Suppose that $\Delta$ is stable. If
$G(I_{\Delta})=\{x_{\sigma_1},x_{\sigma_2},\ldots,x_{\sigma_r}\}$,
then
$$
G(I_{{\em {\rm
Shift}}_{kl}(\Delta)})=\{x_{\sigma_1^*},x_{\sigma_2^*},\ldots,x_{\sigma_r^*}\}.
$$
\end{Lemma}

\begin{pf}
Firstly, we show that $I_{{\rm Shift}_{kl}(\Delta)}$ is generated
by $x_{\sigma_1^*},x_{\sigma_2^*},\ldots,x_{\sigma_r^*}$. It is
enough to show that, for any $\sigma\in\nabla$, there is an $i$
such that $S_{kl}(\sigma)=\sigma^*_i\cup\tau$. As
$G(I_{\Delta})=\{x_{\sigma_1},x_{\sigma_2},\ldots,x_{\sigma_r}\}$,
we have that $\sigma=\sigma_j\cup\tau_1$, for some $j$, where
$\sigma_j<\tau_1$. Divide into three cases.

Case I) $S_{kl}(\sigma_j)\not=\sigma_j$. Then
$\sigma_j^*=S_{kl}(\sigma_j)$, $l\in\sigma_j$, $k\not\in\sigma_j$
and
$S_{kl}(\sigma_j)=(\sigma_j\setminus\{l\})\cup\{k\}\not\in\nabla$,
hence $l\in\sigma$. If $k\in \sigma$, then $k\in\tau_1$ and
\begin{eqnarray*}
S_{kl}(\sigma)
&=&
\sigma\\
&=&
\sigma_j\cup\tau_1\\
&=&
((\sigma_j\setminus\{l\})\cup\{k\})\cup((\tau_1\setminus\{k\})\cup\{l\})\\
&=& \sigma_j^*\cup((\tau_1\setminus\{k\})\cup\{l\}),
\end{eqnarray*}
as required. If $k\not\in\sigma$ and
$(\sigma\setminus\{l\})\cup\{k\}\in\nabla$, then
$(\sigma\setminus\{l\})\cup\{k\}=\sigma_t\cup\tau'$ for some $t$,
where $\sigma_t<\tau'$, and $S_{kl}(\sigma)=\sigma$. Note that
$(\sigma\setminus\{l\})\cup\{k\}=((\sigma_j\setminus\{l\})\cup\{k\})\cup\tau_1$.
Then
$$
((\sigma_j\setminus\{l\})\cup\{k\})\cup\tau_1=\sigma_t\cup\tau',
$$
where $(\sigma_j\setminus\{l\})\cup\{k\}<\tau_1$ and
$\sigma_t<\tau'$. As
$(\sigma_j\setminus\{l\})\cup\{k\}\not\in\nabla$, we see that
$(\sigma_j\setminus\{l\})\cup\{k\}\subseteq \sigma_t$. Thus
$\sigma_t=((\sigma_j\setminus\{l\})\cup\{k\})\cup\tau''=S_{kl}(\sigma_j)\cup\tau''$
with $\tau''\subseteq \tau_1$, then
$\sigma_t^*=\sigma_j\cup\tau''$. Hence
\begin{eqnarray*}
S_{kl}(\sigma) &=&
\sigma\\
&=&
\sigma_j\cup\tau_1\\
&=&
(\sigma_j\cup\tau'')\cup(\tau_1\setminus\tau'')\\
&=& \sigma_t^*\cup(\tau_1\setminus\tau''),
\end{eqnarray*}
as required. If $k\not\in\sigma$ and
$(\sigma\setminus\{l\})\cup\{k\}\not\in\nabla$, then
$S_{kl}(\sigma)=(\sigma\setminus\{l\})\cup\{k\}$. Hence
$S_{kl}(\sigma)=((\sigma_j\setminus\{l\})\cup\{k\})\cup\tau_1=\sigma_j^*\cup\tau_1$,
as required.

Case II) $S_{kl}(\sigma_j)=\sigma_j$ and there is no $s\not=j$
such that $\sigma_j=S_{kl}(\sigma_s)\cup\tau$. Then
$\sigma_j^*=\sigma_j$. If $S_{kl}(\sigma)=\sigma$, then
$S_{kl}(\sigma)=\sigma=\sigma_j^*\cup\tau_1$,as required. If
$S_{kl}(\sigma)\not=\sigma$, then $l\in\sigma$, $k\not\in\sigma$
and $(\sigma\setminus\{l\})\cup\{k\}\not\in\nabla$, hence
$S_{kl}(\sigma)=(\sigma\setminus\{l\})\cup\{k\}=((\sigma_j\cup\tau_1)\setminus\{l\})\cup\{k\}$.
Note that $k\not\in\sigma_j$ implies that $l\not\in\sigma_j$,
otherwise, from $S_{kl}(\sigma_j)=\sigma_j$ we have that
$(\sigma_j\setminus\{l\})\cup\{k\}\in\nabla$, then
$((\sigma_j\cup\tau_1)\setminus\{l\})\cup\{k\}=((\sigma_j\setminus\{l\})\cup\{k\})\cup\tau_1\in\nabla$,
a contradiction. Hence $l\in\tau_1$ and
\begin{eqnarray*}
S_{kl}(\sigma) &=&
\sigma_j\cup((\tau_1\setminus\{l\})\cup\{k\})\\
&=& \sigma_j^*\cup((\tau_1\setminus\{l\})\cup\{k\}),
\end{eqnarray*}
as required.

Case III) $S_{kl}(\sigma_j)=\sigma_j$ and there is some $s\not=j$
such that $\sigma_j=S_{kl}(\sigma_s)\cup\tau$. Then
$\sigma=\sigma_j\cup\tau_1=S_{kl}(\sigma_s)\cup\tau\cup\tau_1$.
Note that $S_{kl}(\sigma_s)\not=\sigma_s$, which implies that
$k\in S_{kl}(\sigma_s)=(\sigma_s\setminus\{l\})\cup\{k\}$. Then
$k\in\sigma$ and $S_{kl}(\sigma)=\sigma$. Since
$\sigma_s^*=S_{kl}(\sigma_s)$, it follows that
$S_{kl}(\sigma)=\sigma=\sigma_s^*\cup(\tau\cup\tau_1)$, as
required.

Secondly, we show that
$x_{\sigma_1^*},x_{\sigma_2^*},\ldots,x_{\sigma_r^*}$ is a minimal
system . It is enough to show that
$\sigma_i^*\not\subseteq\sigma_j^*$ for any $i\not=j$. Divide into
three cases again.

Case 1) $\sigma_j^*=\sigma_j$. It is clear that
$\sigma_i^*\not\subseteq\sigma_j^*$ if $\sigma_i^*=\sigma_i$ or
$\sigma_t\cup\tau$, so we may assume that
$\sigma_i^*=S_{kl}(\sigma_i)$. As $S_{kl}(\sigma_j)=\sigma_j$, if
$\sigma_i^*\subseteq\sigma_j^*$ then $S_{kl}(\sigma_i)\subseteq
S_{kl}(\sigma_j)$, which implies that $\sigma_j^*\not=\sigma_j$, a
contradiction.

Case 2) $\sigma_j^*=S_{kl}(\sigma_j)\not=\sigma_j$. Then, as
$S_{kl}(\sigma_j)\not\in\nabla$, we see that
$\sigma_i^*\not\subseteq\sigma_j^*$ if $\sigma_i^*=\sigma_i$ or
$\sigma_t\cup\tau$. If $\sigma_i^*=S_{kl}(\sigma_i)\not=\sigma_i$
and $\sigma_i^*\subseteq\sigma_j^*$, then $l\in\sigma_i,
\sigma_j$, $k\not\in\sigma_i,\sigma_j$ and
$(\sigma_i\setminus\{l\})\cup\{k\}\subseteq(\sigma_j\setminus\{l\})\cup\{k\}$.
It follows that $\sigma_i\subseteq\sigma_j$, a contradiction.

Case 3) $\sigma_j^*=\sigma_t\cup\tau$, where $\sigma_t<\tau$. Then
$\sigma_j=S_{kl}(\sigma_t)\cup\tau$, $l\in\sigma_t$,
$k\not\in\sigma_t$ and
$S_{kl}(\sigma_t)=(\sigma_t\setminus\{l\})\cup\{k\}\not\in\nabla$.
If $\sigma_i^*=S_{kl}(\sigma_i)\not=\sigma_i$, then
$k\in\sigma_i^*$, hence $\sigma_i^*\not\subseteq\sigma_j^*$ as
$k\not\in\sigma_j^*$. If $\sigma_i^*=\sigma_s\cup\tau'$, where
$\sigma_s<\tau'$, $l\in\sigma_s$ and $k\not\in\sigma_s$, then
$\sigma_i^*\subseteq\sigma_j^*$ implies that
$\sigma_t\cup\tau=\sigma_s\cup\tau'\cup\tau''$ for some $\tau''$,
it follows that
$((\sigma_t\setminus\{l\})\cup\{k\})\cup\tau=((\sigma_s\setminus\{l\})\cup\{k\})\cup\tau'\cup\tau''$,
i.e., $\sigma_j=\sigma_i\cup\tau'\cup\tau''$, a contradiction.
Finally, we assume that $\sigma_i^*=\sigma_i$ and
$\sigma_i^*\subseteq\sigma_j^*$. Then
$\sigma_t\cup\tau=\sigma_i\cup\tau_i$ for some
$\tau_i\not=\emptyset$ and $S_{kl}(\sigma_i)=\sigma_i$. As
$l\in\sigma_t$, $k\not\in\sigma_t$ and $k<l$, we have that
$l\in\sigma_i$ or $l\in\tau_i$, and $k\not\in\sigma_i$ and
$k\not\in\tau_i$. If $l\in\sigma_i$, then, from $k\not\in\sigma_i$
and $S_{kl}(\sigma_i)=\sigma_i$, we have that
$(\sigma_i\setminus\{l\})\cup\{k\}\in\nabla$. Thus
$(\sigma_i\setminus\{l\})\cup\{k\}=\sigma_u\cup\tau_u$ for some
$\sigma_u\in G(I_{\Delta})$. Then, from
$\sigma_t\cup\tau=\sigma_i\cup\tau_i$, we have that
\begin{eqnarray*}
\sigma_j &=& S_{kl}(\sigma_t)\cup\tau\\
&=&
((\sigma_t\setminus\{l\})\cup\{k\})\cup\tau\\
&=& ((\sigma_i\setminus\{l\})\cup\{k\})\cup\tau_i\\
&=&\sigma_u\cup\tau_u\cup\tau_i,
\end{eqnarray*}
a contradiction. If $l\in\tau_i$, then, from
$\sigma_t\cup\tau=\sigma_i\cup\tau_i$, we have that
\begin{eqnarray*}
\sigma_j &=& S_{kl}(\sigma_t)\cup\tau\\
&=& ((\sigma_t\setminus\{l\})\cup\{k\})\cup\tau\\
&=& \sigma_i\cup((\tau_i\setminus\{l\})\cup\{k\}),
\end{eqnarray*}
a contradiction again. The proof is complete.
\end{pf}

Now we can prove

\begin{Theorem}
\label{cbetti} If $\Delta$ is stable, then, for all $i$ and $j$,
$$
\beta_{i,j}(I_{\Delta})=\beta_{i,j}(I_{\Delta^c}).
$$
\end{Theorem}

\begin{pf}
Suppose that
$\Delta^c=\mbox{Shift}_{k_rl_r}(\mbox{Shift}_{k_{r-1}l_{r-1}}(\cdots(\mbox{Shift}_{k_1l_1}(\Delta))\cdots))$
, where $r\geq 1$ and $k_i<l_i$, $i=1,2,\ldots,r$. By
\ref{cstable},
$\mbox{Shift}_{k_il_i}(\mbox{Shift}_{k_{i-1}l_{i-1}}(\cdots(\mbox{Shift}_{k_1l_1}(\Delta))\cdots))$
is stable for $i=1,2,\ldots,r$. Hence it is enough to show that
$\beta_{i,j}(I_{\Delta})=\beta_{i,j}(I_{{\rm
Shift}_{kl}(\Delta)})$ for any $1\leq k<l\leq n$. Let
$G(I_{\Delta})=\{x_{\sigma_1},x_{\sigma_2},\ldots,x_{\sigma_r}\}$,
then, by \ref{genarator}, $ G(I_{{\rm
Shift}_{kl}(\Delta)})=\{x_{\sigma_1^*},x_{\sigma_2^*},\ldots,x_{\sigma_r^*}\}$.
Since $I_{\Delta}$ and $I_{{\rm Shift}_{kl}(\Delta)}$ are
squarefree stable and $m(\sigma_i^*)=m(\sigma_i)$,
$i=1,2,\ldots,r$, it follows from the formula of Betti numbers of
squarefree stable ideals that
$\beta_{i,j}(I_{\Delta})=\beta_{i,j}(I_{{\rm
Shift}_{kl}(\Delta)})$, as required.
\end{pf}

Combining \ref{sbetti}, \ref{ebetti} and \ref{cbetti}, we get the
following

\begin{Theorem}
\label{betti} If $\Delta$ is stable, then, for all $i$ and $j$,
$$
\beta_{i,j}(I_{\Delta^s})=\beta_{i,j}(I_{\Delta^e})=\beta_{i,j}(I_{\Delta^c}).
$$
\end{Theorem}

\medskip
\medskip
\noindent
Zhongming Tang\\
Department of Mathematics\\
Suzhou University\\
Suzhou 215006\\
P.\ R.\ China\\
E-mail: zmtang@@suda.edu.cn

\medskip
\noindent
Guifen Zhuang\\
Department of Mathematics\\
Suzhou University\\
Suzhou 215006\\
P.\ R.\ China\\




\end{document}